\documentclass[12pt]{amsart}
\usepackage{amscd,amssymb,url}
\usepackage[graph,frame,poly,arc]{xy}  
\usepackage[plainpages,backref,urlcolor=blue]{hyperref}

\usepackage{tikz}
\usetikzlibrary{calc}

\theoremstyle{plain}
\newtheorem{thm}[subsection]{Theorem}

\newtheorem{prop}[subsection]{Proposition}
\newtheorem{cor}[subsection]{Corollary}

\theoremstyle{definition}
\newtheorem{rk}[subsection]{Remark}
\newtheorem{definition}[subsection]{Definition}

\newtheorem{conj}[subsection]{Conjecture}

\numberwithin{equation}{section}
\newcommand{\HH}{{\mathcal H}}

\newcommand{\F}{{\mathcal F}}
\newcommand{\A}{{\mathcal A}}
\newcommand{\B}{{\mathcal B}}

\newcommand{\M}{{\mathcal M}}

\newcommand{\CC}{{\mathcal C}}
\newcommand{\LL}{{\mathcal L}}

\newcommand{\al}{{\alpha}}
\newcommand{\be}{{\beta}}
\newcommand{\la}{{\lambda}}

\newcommand{\PPP}{{\mathcal P}}

\newcommand{\C}{\mathbb{C}}
\newcommand{\FF}{\mathbb{F}}
\newcommand{\PP}{\mathbb{P}}

\begin{document}

\title [On the Alexander polynomials of conic-line arrangements]
{On the Alexander polynomials of conic-line arrangements}

\author[Alexandru Dimca]{Alexandru Dimca}
\address{Universit\'e C\^ ote d'Azur, CNRS, LJAD, France and Simion Stoilow Institute of Mathematics,
P.O. Box 1-764, RO-014700 Bucharest, Romania}
\email{Alexandru.Dimca@univ-cotedazur.fr}

\author[Piotr Pokora]{Piotr Pokora}
\address{Department of Mathematics,
University of the National Education Commission Krakow,
Podchor\c a\.zych 2,
PL-30-084 Krak\'ow, Poland.}
\email{piotr.pokora@uken.krakow.pl, piotrpkr@gmail.com}

\author[Gabriel Sticlaru]{Gabriel Sticlaru}
\address{Faculty of Mathematics and Informatics,
Ovidius University
Bd. Mamaia 124, 900527 Constanta,
Romania}
\email{gabriel.sticlaru@gmail.com}

\thanks{\vskip0\baselineskip
\vskip-\baselineskip
\noindent A. Dimca was partially supported by the Romanian Ministry of Research and Innovation, CNCS - UEFISCDI, Grant \textbf{PN-III-P4-ID-PCE-2020-0029}, within PNCDI III.\\
P. Pokora is partially supported by The Excellent Small Working Groups Programme \textbf{DNWZ.711/IDUB/ESWG/2023/01/00002} at the Pedagogical University of Cracow. }

\subjclass[2010]{Primary: 14H50; Secondary:  14N20, 14B05, 13D02, 32S35, 32S40, 32S55}

\keywords{plane curves; conic pencils; free curves;  Alexander polynomials}

\begin{abstract} 
In the present paper we compute Alexander polynomials for certain classes of conic-line arrangements in the complex projective plane
which are related to pencils. We prove two general results for curve arrangements coming from Halphen pencils of index $k\geq 2$. 
Then we apply them to the Hesse arrangement of conics and to some of its degenerations. The results are completed by computations using computer algebra. In particular, we construct conic-line arrangements which are non-reduced pencil-type arrangements and have as roots of their Alexander polynomials roots of unity of order 7. Such roots are not known and are conjectured not to exist in the class of line arrangements.
\end{abstract}
 
\maketitle

\section{Introduction} 

Let $S=\C[x,y,z]$ be the graded polynomial ring in the variables $x,y,z$ with the complex coefficients and let $C : f=0$ be a reduced curve of degree $d$ in the complex projective plane $\mathbb{P}^{2}=\mathbb{P}^2_{\mathbb{C}}$. 
If the curve $C$ is reducible, it is called a curve arrangement. 
Let $g_1,g_2 \in S_m$ be two  homogeneous polynomials of degree $m$, without common factor, and consider the associated pencil of degree $m$ curves
$$\PPP:sg_1+tg_2=0,$$
where $(s:t) \in \PP^1$. We assume that the generic member of $\PPP$ is irreducible.
Such a pencil is called a {\it Halphen pencil of index } $k$ if one has $g_2=h^k$ for an integer $k>1$ and a polynomial $h$ which is not a power of another polynomial.
A curve $C$ is called a {\it reduced (resp. non-reduced) pencil-type arrangement associated with $\ell$ members} of $\PPP$ if there is a number $\ell\geq 3$  of reduced members (resp. members with at least one member having a multiple irreducible component) of the pencil $\PPP$ such that $C$ is the union of these $\ell$ members (after replacing all the multiple components by their reductions). For instance, the monomial arrangement
$$\M_m: (x^m-y^m)(y^m-z^m)(x^m-z^m)=0$$
is a reduced pencil-type line arrangement with $\ell=3$ members associated to the pencil
$\PPP: s(x^m-y^m)+t(y^m-z^m)=0$, while the full monomial arrangement
$$\F\M_m:xyz(x^m-y^m)(y^m-z^m)(x^m-z^m)=0$$
is a non-reduced pencil-type line arrangement with $\ell=3$ associated to the pencil
$\PPP: sz^m(x^m-y^m)+tx^m(y^m-z^m)=0$. We do not know if a curve $C$ can be a reduced  pencil-type arrangement and, in the same time,
a  non-reduced pencil-type arrangement, possibly with respect to two distinct pencils. However, this fact does not affect our results below.

Such pencil-type arrangements play a key role in the theory of line arrangements, see for instance \cite{FY, DHA}, and their relation to the freeness was considered in \cite{DMich, JV}.

We consider the complement  $U= \PP^2 \setminus C$ and let $F:f=1$ be the corresponding Milnor fiber in $\C^3$, with the usual monodromy action $h:F \to F$. One can also consider  the characteristic polynomials of the monodromy, namely
\begin{equation} 
\label{Delta}
\Delta^j_C(t)=\det (t\cdot Id -h^j \,|\, H^j(F,\C)),
\end{equation} 
for $j \in \{0,1,2\}$. Since the curve $C$ is reduced,  one has $\Delta^0_C(t)=t-1$, and moreover
\begin{equation} 
\label{Euler}
\Delta^0_C(t)\Delta^1_C(t)^{-1}\Delta^2_C(t)=(t^d-1)^{\chi(U)},
\end{equation} 
where $\chi(U)$ denotes the Euler characteristic of the complement $U$, see for instance  \cite[Proposition 4.1.21]{D1}. Let us recall that
\begin{equation} 
\label{Euler1}
 \chi(U)=(d-1)(d-2)+1-\mu(C),
\end{equation}  
where $\mu(C)$ is the total Milnor number of $C$, which is the sum of the Milnor numbers $\mu(C,p)$ for all the singular points $p$ of $C$. 
 It follows that the polynomial $\Delta(t)=\Delta^1_C(t)$, also called the Alexander polynomial of $C$ -- see \cite{Li1,Li6,Ra}, determines the remaining polynomial $\Delta^2_C(t)$ if one knows $\mu(C)$. 
 
 When $\chi(U) \leq 0$, it is known that all the irreducible components $C_i$ of $C$ are rational,  see \cite{GuPa, JoSt, WV}, and then $\Delta(t)$ tends to be quite large. If these components $C_i$ are smooth,
 being rational means that they are either lines or conics, so $C$ is a conic-line arrangement.
 However, it is easy to show that for any line arrangement $C$ with $d=\deg(C) \geq 5$ and without points of multiplicity $\geq d-1$ one has
 \begin{equation} 
\label{EulerA}
 \chi(U) \geq d-4 \geq 1,
\end{equation} 
 see Proposition \ref{propA}. 
 Note that the arrangements of $d$ lines having points of multiplicity $\geq d-1$ or such that $d<5$ are easy to describe, and their monodromy is well-known. In fact, for most line arrangements $C$, the Alexander polynomial
 $\Delta^1_C(t)$ is trivial, namely it satisfies
\begin{equation} 
\label{triv1}
\Delta^1_C(t)=(t-1)^{r-1}, 
\end{equation}  
where $r$ is the number of irreducible components of $C$, see for instance \cite{CS,MP,S}. Even when the Alexander polynomial
 $\Delta^1_C(t)$ is non-trivial, possible roots of  $\Delta^1_C(t)$ seem to be very restricted. For instance, we have the following conjecture by Papadima and Suciu \cite[Conjecture 1.9]{PS}.
\begin{conj} 
\label{conj1} 
Let $\al$ be a root of the Alexander polynomial 
 $\Delta^1_C(t)$ of a line arrangement $C$, not a union of concurrent lines.
 If the order of $\al$ has the form $p^s$, with $p$ a prime number and
 $s>0$ an integer, then $p^s \in \{2,3,4\}$.
\end{conj} 
Moreover, all known examples of line arrangements with non trivial Alexander polynomials come from pencil-type arrangements, see \cite{FY,S,S4}. More precisely, if a {\it reduced pencil-type line arrangement} consists of $\ell$  members of a pencil, then the corresponding Alexander polynomial has as roots any primitive root of order $\ell$, see Corollary \ref{corPEN}. 
On the other hand, for a {\it non-reduced pencil-type line arrangement}
consisting of $\ell$ members, the roots of unity of order $\ell$ are not always roots of its Alexander polynomial. For instance, for the full monomial arrangement $\F\M_m$ the cubic roots of unity are roots of the Alexander polynomial if and only if $m \equiv 1$ mod $3$, see \cite[Remark 3.4]{DP}.
Since for any line arrangement of pencil-type, except the union of concurrent lines, one has $\ell \leq 4$, this was perhaps the motivation for Conjecture \ref{conj1}.
The first author has constructed in \cite{JKMS} some conic-line arrangements $C$, related to Halphen pencils of conics of index 2, with $\chi(U) \leq 1$ and such that 
$\Delta^1_C(t)$ admits as roots any root of unity of order $d=\deg(C)$, see Remark \ref{rkPENCILS}  for additional information.

In this note we  construct (reduced and non-reduced) pencil-type conic-line arrangements related to pencils of conics and sextics, with non-trivial Alexander polynomial
 $\Delta^1_C(t)$ and such that $\chi(U)>0$. For instance, the Hesse arrangement of conics $C=\mathcal{C}_{12}(\lambda)$ from Theorem \ref{thmHesse} satisfies $\chi(U)=54$, consists of $\ell=4$ reduced members in a pencil of sextics, and $\Delta^1_C(t)$ has roots of order $8$. This arrangement has non quasi-homogeneous singularities, but which are locally topological equivalent to quasi-homogeneous singularities.
It seems therefore that Conjecture  \ref{conj1} has nothing to do with the positivity of $\chi(U)$ or the local topological properties of singularities, but is likely to be a property specific to line arrangements.
 
It is an interesting question to see how the freeness of a curve $C$ is reflected in the topological properties of $U$ and $F$, see for instance \cite{ArD}.
In this note we show, in particular, that some new free and nearly free conic-line arrangements lead to arrangements with non-trivial Alexander polynomial
 $\Delta^1_C(t)$.
The freeness of conic-line arrangements is also discussed in \cite{STo}, from a different perspective and with an aim towards Terao's freeness conjecture.

 Let us briefly present an outline of our note. In Section 2, we review some basic facts on the Alexander polynomials of curve arrangements $C$ and prove the inequality \eqref{EulerA} in Proposition \ref{propA}.
 In Section 3, we discuss the relation between Alexander polynomials and pencil-type arrangements. We recall the basic facts in Theorem \ref{thmPEN} and prove two of our main results, namely Theorem
 \ref{thmPEN2} and Theorem \ref{thmPEN3} which explain the additional information provided by Halphen pencils. Using these results, we can predict the {\it existence of certain roots of the Alexander polynomials in all the cases discussed in this paper}, see Remarks \ref{rkPENCILS}, \ref{rkHesse2P} and \ref{rkHesse3P}.
 
 In Section 4, we consider conic-line arrangements associated to pencils of conics. Their freeness, except from the case of the generic pencil of Type I, see Theorem \ref{thm1}, is not a consequence of the results in \cite{DMich, JV}.
Both the freeness and the Alexander polynomials are computed using
the algebra software \verb}Singular} via the algorithms developed in
\cite{DStFor,DStMFgen,DStproj}, which are recalled in Section 2.

In Section 5, we discuss the main example of our paper, the $1$-parameter family of Hesse arrangements of conics\footnote{Sometimes this arrangement is called as the Chilean arrangement of conics.} $\mathcal{C}_{12}(\lambda)$ for 
$$\lambda \in \C \setminus \{-2,-2 \omega, -2 \omega^2, 0, 1, \omega,  \omega^2 \}$$
with $ \omega^2+ \omega+1=0$, of $12$ smooth conics in the plane \cite{Dolgachev, KRS}, which is both free and has an interesting Alexander polynomial.  In Section 6, we consider some degenerations of the Hesse arrangement of conics $\mathcal{C}_{12}(\lambda)$, obtained by setting $\lambda=\lambda_0$ for 
$$\lambda_0 \in  \{-2,-2 \omega, -2 \omega^2, 0, 1, \omega,  \omega^2 \}$$
in the equation $f(\lambda)=0$ of the aforementioned $1$-parameter family of Hesse arrangements of conics.
For each of these values of $\lambda$, the defining equation $f(\lambda)=0$ has 3 linear factors which occurs with multiplicity 2, hence $f(\lambda)=0$ is in fact a multiarrangement $\CC(\lambda)$ formed by $9$ smooth conics and $3$ double lines -- here by a multiarrangement we understand an arrangement of curves such that not every irreducible component is reduced. From the point of view of pencils of sextics, these degenerations  are union of 4 members (one double and 3 reduced) of a Halphen pencil of sextics of index 2.
We denote  by $\A(\lambda)$ the {\it non-reduced pencil-type} conic-line arrangement
with is obtained from $\CC(\lambda)$ by replacing the 3 double lines by the corresponding reduced lines, and by $\B(\lambda)$ the {\it reduced pencil-type} subarrangement of $\A(\lambda)$ consisting just of the 3 reduced members.
These 14 new arrangements $\A(\lambda)$ and $\B(\lambda)$ have very interesting properties themselves. For instance, Theorem \ref{thmHesse1} says that all the seven conic-line arrangements $\A(\lambda)$ are free, and only three of the seven conic arrangements
$\B(\lambda)$, namely those corresponding to $\lambda^3=1$, are free.
Finally, Theorem \ref{thmHesse2} and Theorem \ref{thmHesse3} say that unit roots of order 6 (resp. 7) occur as roots of the Alexander polynomial of the arrangements $\B(0)$ and $\B(1)$ (resp. $\A(0)$ and $\A(1)$).
It is an interesting fact that the combinatorics and the local type of singularities of the arrangements
 $\A(0)$ and $\A(1)$ (resp. $\B(0)$ and $\B(1)$) are very different, but nevertheless one has
 $$\Delta^1_{\A(0)}(t)= \Delta^1_{\A(1)}(t) \text{ and } \Delta^1_{\B(0)}(t)= \Delta^1_{\B(1)}(t).$$
One explanation of this fact comes by looking at the associated pencils of sextics: both $\B(0)$ and $\B(1)$ consist
of $3$ reduced fibers, and $\A(0)$ and $\A(1)$ are obtained by adding the  reduce curve coming from the double cubic in the pencil, which is in fact a triangle.

We would like to thank the referees for their careful reading of the manuscript and for their very useful suggestions to improve the presentation.

 \section{Some general properties of Alexander polynomials of curve arrangements} 

In order  to construct conic-line arrangements $C$ with non-trivial Alexander polynomial
 $\Delta^1_C(t)$, the following result is quite useful, see \cite[Section 6.4]{Sheaves}. We recall that for any isolated plane curve singularity $(X,0)$, there is a Milnor fiber $F_0$ and a monodromy operator $h_0:F_0 \to F_0$, giving rise to a local Alexander polynomial
 \begin{equation} 
\label{Delta0}
\Delta^1_{(X,0)}(t)=\det (t\cdot {\rm Id} - h^1_0 \,|\, H^1(F_0,\C)).
\end{equation}

 \begin{thm}
\label{thm0} 
Let $C=\cup_{i=1}^sC_i$ be a curve arrangement in $\PP^2$ such that the curve $C_1$ is irreducible. For any singular point $p$ of $C$ situated on $C_1$, let $ \Delta^1_{(C,p)}(t)$ be the  Alexander polynomial of the singularity $(C,p)$. Then the following conditions hold.
\begin{enumerate}

\item  Any root of the Alexander polynomial
 $\Delta^1_C(t)$ is a root of unity of order $d=\deg(C)$.

\item  The Alexander polynomial
 $\Delta^1_C(t)$ divides the product $\prod_p  \Delta^1_{(C,p)}(t)$, where
 $p$ runs through all the singularities of $C$ situated on $C_1$.

\end{enumerate}

\end{thm}
The claim (i) above is classical, going back to Zariski. For the claim (ii) we refer to \cite[Corollary 6.4.16]{Sheaves}.

\begin{rk}
\label{rk0}  
Here are some examples of local Alexander polynomials. For the simple singularity $A_n:u^2+v^{n+1}=0$, the corresponding Alexander polynomial is 
$$\Delta^1_{(A_n)}(t)=\frac{t^{n+1}-(-1)^{n+1}}{t+1}.$$
For an ordinary $m$-multiple point $(X_m,0)$, that is when there are $m$ smooth branches meeting transversally in pairs at a given point, we have
$$\Delta^1_{(X_m,0)}(t)=(t^{m}-1)^{m-2}(t-1).$$
To get these formulas, we refer to the proof of Theorem (3.4.10) in \cite{D1}.
\end{rk}


Theorem \ref{thm0} gives necessary conditions to get interesting eigenvalues for the monodromy. However, these conditions are far from being sufficient. To make sure that a root of unity $\al$ such that $\al^d=1$ is a root of the Alexander polynomial $\Delta^1_C(t)$, we can follow {\it two approaches}. 

{\it The first approach} is to compute the Alexander polynomials of plane curves  using two Singular codes, which can be found at 
 \begin{center}
\url{https://math.unice.fr/~dimca/singular.html}
\end{center}

The first one, the code \verb}monof3(q1,q2)}, computes the Alexander polynomial  $\Delta^1_C(t)$  for a {\it free} curve $C$ and is very rapid. Its construction is described in  \cite{DStFor, DHA}, see also \cite[Remark 4.4]{DStMFgen}. 

The second one, the code \verb}mono3_wh()}, described in  \cite{DStMFgen, DHA}, computes the Alexander polynomial  $\Delta^1_C(t)$  for a {\it any} reduced curve $C$, i.e., we do not require the freeness, etc. Moreover, it computes the total Milnor number $\mu(C)$ and then, using the equalities \eqref{Euler} and \eqref{Euler1}, one can determine the Alexander polynomial  $\Delta^2_C(t)$ as well (which is in fact also given by this code). However, for this second code, make sure that the coordinates $x,y,z$ are chosen such that the line at infinity $L:z=0$ contains no singularities of the curve $C$. This is achieved by making a linear substitution $z \mapsto ax+by+z$ with $a,b$ being general enough.
Otherwise the computation of $\mu(C)$ gives a wrong value.

More precisely, let 
 \begin{equation} \label{root}
\al_q=\exp(-2 \pi iq/d)
\end{equation} 
 be a root of unity of order $d= \deg C$, where $0 \leq q \leq d$ and let $m(\al_q)$ be the multiplicity of
$\al_q$ as a root of the Alexander polynomial $\Delta^1_C(t)$. 
One has $\al_0=\al_d=1$ and 
 \begin{equation} \label{mult0}
m(1)= b_1(U)=r-1,
 \end{equation} 
where $r$ is the number of irreducible components of $C$ and $b_{1}(U)$ denotes the first Betti number of the complement $U$ of $C$.
The output of the code \verb}monof3(q1,q2)} with $q1=3$ and $q2=d$ is a Table, where the first column contains the values of $q$ from $3$ to $d$,
and the before last column contains the data
$n_2(q)$. We set $n_2(0)=n_2(1)=n_2(2)=0$. Then \cite[Remark 4.4]{DStMFgen} yields
 \begin{equation} \label{mult}
m(\al_q)=n_2(q)+n_2(d-q) \text{ for any } q=0, \ldots, d.
\end{equation} 

The output of the code \verb}mono3_wh()}, when the curve has some non quasi-homogeneous singularities, that is when $\tau(C) < \mu(C)$, is first the total Tjurina number $\tau(C)$ and the total Milnor $\mu(C)$. Then  there is again a Table, where the first column contains the values of $q$ from $3$ to $d$,
and the before last column contains the data
$H^1_q$. We set $H^1_0=H^1_1=H^1_2=0$. Then \cite[Remark 4.4]{DStMFgen} yields
 \begin{equation} 
\label{mult2}
m(\al_q)=H^1_q+H^1_{d-q} \text{ for any } q=0, \ldots, d.
\end{equation} 
When all the singularities of $C$ are quasi-homogeneous, then there are additional information in the output, but for the computation of  $m(\al_q)$ one can use the data displayed in the first part of Table that is produced by the code \verb}mono3_wh()} and the formula \eqref{mult2}. In view of \cite[Proposition 2.2 (2)]{DStproj} and \cite[Remark 8.6]{DHA}, the numbers
$n_2(q)$ and $H^1_q$ have the following geometric interpretation
 \begin{equation} 
\label{mult20}
n_2(q)= \dim Gr^1_FH^1(F,\C)_{\al_q},
\end{equation} 
when all the singularities of $C$ are quasi-homogeneous, and
 \begin{equation} 
\label{mult21}
H^1_q= \dim Gr^1_FH^1(F,\C)_{\al_q},
\end{equation}
otherwise, where $Gr^1_F$ denotes the graded piece with respect to the Hodge filtration on $H^1(F,\C)$.

The interested reader can find alternative codes for the computation of Alexander polynomials obtained by Morihiko Saito in \cite{Sa3,Sa4}.

\medskip

{\it The second approach} is to use the following formula for the multiplicity $m(\al)$ of $\al$ as a root of the Alexander polynomial $\Delta_C(t)$, namely
 \begin{equation} 
\label{mult3}
m(\al)=\dim H^1(U,\LL_{\al}),
\end{equation} 
where $\LL_{\al}$ is the rank one local system whose monodromy around each irreducible component of $C$ is multiplication by $\al$,
see \cite[Proposition 5.4]{DHA}. This approach is particularly useful when
the curve $C$ is the union of some members in a pencil of curves, as we explain in the next section.

We end this section with the following easy result.

 \begin{prop}
\label{propA} 
Let $C:f=0$ be an arrangement of $d$ lines in $\PP^2$, without points of multiplicity $\geq d-1$. If $d=\deg(f) \geq 5$ and $U= \PP^2 \setminus C$, then 
$$\chi(U) \geq d-4 \geq 1.$$
 \end{prop}
 \proof
Since $C$ has no points of multiplicity $\geq d-1$, one has $r={\rm mdr}(f) \geq 2$ (regarding $r$, see discussion below Definition \ref{m-syz}).
Furthermore, the singular points of $C$ are quasi-homogeneous, so we have
$$\mu(C)=\tau(C)\leq \tau(d,r)_{max} \leq \tau(d,2)_{max}=d^2-4d+7,$$
and this follows from \cite{duP}. Hence we have
$$\chi(U)=(d-1)(d-2)+1-\mu(C)\geq (d-1)(d-2)+1-(d^2-4d+7)=d-4.$$
 \endproof
 Note that the arrangements of $d$-lines having points of multiplicity $\geq d-1$ or such that $d<5$ are easily described, see for instance \cite[Cor. 8.1 and Cor. 8.2]{DHA}, and their monodromy is well-known.

\section{Pencils of curves and Alexander polynomials} 
Let $g_1,g_2 \in S_m$ be two linearly independent homogeneous polynomials of degree $m$, consider the associated pencil
$\PPP:sg_1+tg_2=0,$
and assume that the generic member of $\PPP$ is irreducible.
If the curve $C$ is a {\it pencil-type arrangement associated with $\ell\geq 3$  (reduced or non-reduced) members} of $\PPP$, then  we have a regular mapping
 \begin{equation} 
\label{p1}
\theta:U \to T, \  \ \theta(x,y,z)=(g_1(x,y,z):g_2(x,y,z)),
\end{equation}
where $T$ is obtained from $\PP^1$ by deleting the $\ell$ points corresponding to the $\ell$  members of $\PPP$ used to construct $C$. Note that one has 
$ \chi(T)=2- \ell.$
We recall the following key result, see \cite[Theorem 6.9, Corollary 6.3 and Theorem 6.10]{DHA}. For all the details we refer to \cite{D4}.
\begin{thm}
\label{thmPEN} 
Let $\LL$ be a rank one local system on $U$ and let $\FF$ be the generic fiber of $\theta$. If $\LL|_{\FF}=\C_{\FF}$ is the trivial local system, 
then for any rank one local system $\LL'$ on $T$ one has
$$\dim H^1(U,\LL \otimes \theta^*\LL') \geq -\chi(T)=\ell -2 \geq 1.$$
The first inequality is strict when $\LL$ is not the pull-back under $\theta$ of a rank one local system on $T$.
\end{thm}
\begin{cor}
\label{corPEN} 
Any root of unity $\al$ of order $\ell$ is a root for the Alexander polynomial of a  reduced pencil-type arrangement associated with $\ell$ members of $\PPP$ and moreover $m(\al) \geq \ell-2$.
\end{cor}
\proof
It is enough to take in Theorem \ref{thmPEN} $\LL=\C_U$ and $\LL'$ the rank one local system on $T$ with monodromy about each deleted point a primitive root of unity of order $\ell$.
\endproof
When $\PPP$ is a {Halphen pencil of index} $k$, then the map $\theta$ has a multiple fiber, namely $F_0=\theta^{-1}(1:0)$, of order $k$. Then the following result shows the existence of roots of Alexander polynomial of higher order and higher multiplicity for the associated reduced pencil-type arrangements.
\begin{thm}
\label{thmPEN2} 
Let $\PPP$ be a Halphen pencil of index $k$ and let $C$ be a reduced pencil-type arrangement associated with $\ell$  members of $\PPP$. Then the Alexander polynomial $\Delta_C(t)$ has, among its roots, any  primitive root of unity $\al$ of order $k\ell$ and its multiplicity satisfies $m(\al)\geq \ell-1.$
\end{thm}
\proof
We apply Theorem \ref{thmPEN} in the following way. Let $\al$ be any primitive root of unity of order $k\ell$ and consider the rank one local system $\LL_0'$ on $T'=T\setminus \{(1:0)\}$ whose monodromies are as follows:
\begin{enumerate}

\item About each of the $\ell$ points corresponding to the member of
$\PPP$ which are in $C$ the monodromy is multiplication by $\al$.

\item  About the point $(1:0)$ the monodromy is $\be=\al^{-\ell}$.

\end{enumerate}
Since the product of all these monodromies is equal to $1$, this local system $\LL_0'$ is well-defined. Let $\LL_1'=\theta^*(\LL_0')$ be the pull-back of $\LL_0'$ under the obvious restriction of $\theta$ to $U'=\theta^{-1}(T')$. Each irreducible component $C'$ of the multiple fiber
$F_0$ has a multiplicity of the form $ka$ for some integer $a\geq 1$.
It follows that the monodromy of the local system $\LL_1'$ about $C'$ is
$$\be ^{ka}=\al^{-\ell ka}=1.$$
It follows that the local system $\LL_1'$ can be extended over the fiber $F_0$, that is, there is a local system $\LL_{\al}$ on $U$ whose restriction to $U'$ is exactly $\LL_1'$. It is clear that this local system $\LL_{\al}$ satisfies the equality \eqref{mult3}. On the other hand, $\LL_{\al}$ is obviously not the pull-back of a rank one local system on $T$. The claims then follow from Theorem \ref{thmPEN}, where we take $\LL=\LL_{\al}$ and $\LL'=\C_T$. The property $\LL|_{\FF}=\C_{\FF}$ holds since the generic fiber of $\theta$ is the same as the generic fiber of $\theta'=\theta |_{U'}$.
 
\endproof
 We also have the following result on the simplest non-reduced pencil-type arrangements.
 
 \begin{thm}
\label{thmPEN3} 
Let $\PPP$ be a Halphen pencil of index $k$ given by
$\PPP:sg_1+tg_2=0,$ with $g_2=h^k$ and $h$ a polynomial without multiple factors.
Let $C$ be a non-reduced pencil-type arrangement associated with $\ell-1\geq 2$  reduced members of $\PPP$ and the reduced curve $h=0$. Then the Alexander polynomial $\Delta_C(t)$ has, among its roots, any root of unity $\al$ of order $n=k(\ell-1)+1$ and its multiplicity satisfies $m(\al)\geq \ell-2.$
\end{thm}
\proof
Let $\al$ be any  root of unity  of order $n$ and consider the rank one local system $\LL'$ on $T$ whose monodromies are as follows:
\begin{enumerate}

\item  About each of the $\ell-1$ points corresponding to the  reduced members of
$\PPP$ which are in $C$ the monodromy is multiplication by $\al$.

\item  About the point $(1:0)$, corresponding to the multiple member, the monodromy is multiplication by $\be=\al^{1-\ell}$.

\end{enumerate}
We have  $\be^k=\al^{k-k\ell}=\al$, which implies that $\theta^*(\LL')=\LL_{\al}$, with the notation from \eqref{mult3}.
We conclude using Theorem \ref{thmPEN} for $\LL=\C_U$, the trivial rank one local system. This gives in particular 
$m(\al) \geq -\chi(T)=\ell-2$.
\endproof

\begin{rk}
\label{rkGF}  
If the map $\theta: U \to T$ associated to the pencil $\PPP:sg_1+tg_2=0$ as in \eqref{p1} has a generic fiber $\FF$ which is not connected, then Stein Factorization Theorem implies that the rational map $\hat \theta : \PP^2 \dasharrow \PP^1$ associated to $\PPP$ factors as a composition $\PP^2 \dasharrow \PP^1 \to \PP^1$, where
the rational map $\PP^2 \dasharrow \PP^1$ comes from a lower-degree pencil $\PPP'$ with a connected generic fiber $\FF'$ and $\eta:\PP^1 \to \PP^1$ is a finite morphism. Since any member of the pencil $\PPP$ is a union of members of the pencil $\PPP'$, it follows that any arrangement associated with $\PPP$ can be regarded as an arrangement associated to $\PPP'$. Hence our condition that $\FF$ is connected is not really a restriction. To have an example, if the pencil $\PPP$ has at least two fibers of multiplicity $k$, then one can take $g_1=h_1^k$ and $g_2=h_2^k$. In this case, if say the curve $h_1=0$ is irreducible, then the  pencil $\PPP'$ is given by
$sh_1+th_2=0$  and 
$\eta:\PP^1 \to \PP^1$ is given by $\eta(u:v)=(u^k:v^k)$.
\end{rk}

\section{Free curves and pencils of conics} 

In this section we recall first the definition of free and nearly free curves.
In the sequel, many of the constructed curves have not only interesting Alexander polynomials, but turn out to be either free or nearly free.
 Let $C : f=0$ be a reduced curve in $\mathbb{P}^{2}$ of degree $d$ defined by $f \in S$. Denote by $M(f) := S/ J_{f}$ the Milnor algebra, where $J_{f}$ is the Jacobian ideal, i.e., the ideal spanned by the partials $\partial_{x}f, \partial_{y}f, \partial_{z}f$. 
\begin{definition}
\label{m-syz}
We say that $C$ is $m$-syzygy when $M(f)$ has the following minimal graded free resolution
$$0 \rightarrow \bigoplus_{i=1}^{m-2}S(-e_{i}) \rightarrow \bigoplus_{i=1}^{m}S(1-d - d_{i}) \rightarrow S^{3}(1-d)\rightarrow S \rightarrow M(f) \rightarrow 0$$
with $e_{1} \leq e_{2} \leq \ldots \leq e_{m-2}$ and $1\leq d_{1} \leq \ldots \leq d_{m}$.
\end{definition}
In the setting of the above definition, we define the minimal degree of the Jacobian relations among the partial derivatives of $f$, namely
$${\rm mdr}(f) := d_{1}.$$
\begin{definition}
We say that 
\begin{itemize}
\item $C$ is free if and only if $m=2$ and $d_{1}+d_{2}=d-1$. 
\item $C$ is nearly-free if and only if $m=3$, $d_{1}+d_{2} = d$, $d_{2}=d_{3}$. 
\end{itemize}
In both cases, we call the pair $(d_1,d_2)$ the exponents of the curve $C$.
Finally, for a reduced curve $C\subset \mathbb{P}^{2}_{\mathbb{C}}$ let us denote by $\tau(C)$ the total Tjurina number of $C$ and by $\mu(C)$ the total Milnor number of $C$. If all the singularities of $C$ are quasi-homogeneous, then, according to \cite{Reiffen}, one has $\tau(C) = \mu(C)$.
There are characterizations of free and nearly free curves in terms of
the total Tjurina number -- see \cite{Dmax,duP} for all necessary details.

In  the remaining of this section we apply Theorem \ref{thmPEN} to curves coming from pencils of conics.
We start by recalling the classification of pencils of conics $\PPP: tq_1+sq_2=0$ in $\PP^2$ which contain smooth conics. Let $B:q_1=q_2=0$ be the base locus, regarded as a scheme in $\PP^2$, and let $s$ and $s'$ be the number of conics which are union of two distinct lines in $\PPP$, respectively the number of conics which are double lines in $\PPP$. We will use a rather standard convention of describing base loci, for instance by $B = p_{1} + p_{2} + p_{3} + p_{4}$ we mean the base locus $B$ consisting of $4$ distinct points. Then the following result is well-known, see \cite[Section 3]{AEI} or \cite{Per}.

\begin{thm}
\label{thm1} 
With the above notation, the following $5$ cases occur for a pencil of conics $\PPP$ in $\PP^2$ which contains smooth conics.
\begin{enumerate}

\item  Type I (the generic case): $B=p_1+p_2+p_3+p_4$, $s=3$, $s'=0$.

\item  Type II: $B=2p_1+p_2+p_3$, $s=2$, $s'=0$.

\item Type III: $B=3p_1+p_2$, $s=1$, $s'=0$.

\item Type IV: $B=2p_1+2p_2$, $s=s'=1$.

\item Type V: $B=4p_1$, $s=0$, $s'=1$.

\end{enumerate}

\end{thm}

\begin{rk}
\label{rk1}  Here are some equations for each type of pencil of conics listed above, see \cite{Per} for all types except Type I.
$${\rm Type \, I}: q_1=x^2-y^2, q_2=y^2-z^2,$$
and the third singular conic in this pencil is $q_3=x^2-z^2$.
$${\rm Type\, II}:  q_1=xy, q_2=(x-y)(x+y+z),$$
and the common tangent at the point $p_1=(0:0:1)$ is $T:x-y=0$.
$${\rm Type\, III}:  q_1=x(x-z), q_2=xy-z^2.$$
$${\rm Type\, IV}:  q_1=(x+y+z)^2, q_2=xy,$$
and the common tangents to the conics are $T_1:x=0$ and $T_2:y=0$.
$${\rm Type\, V}:  q_1=x^2, q_2=xy-z^2.$$
Note that a pencil of Type IV is called a {\it bitangent pencil} and a pencil of Type V is called a {\it hyperosculating pencil}, see for instance \cite{Be}. The monodromy of some conic-line arrangements constructed using such pencils of Type IV and V has already been discussed in \cite{JKMS}.
\end{rk}

Now we are ready to present our results devoted to pencils of conics and their Alexander polynomials.
\end{definition}
\begin{prop}
\label{propT1a} 
Consider the conic-line arrangement 
$$C:f=(x^2-y^2)(y^2-z^2)(x^2-z^2)(x^2-2y^2+z^2)=0$$ 
obtained from a conic pencil of Type I by putting together the 6 lines coming from the 3 singular conics in the pencil, plus a smooth member. Then $C$ is a free curve with exponents $(2,5)$ and one has 
$\tau(C)=\mu(C)=39$ and 
$$\Delta^1_C(t)=(t-1)^4(t^4-1)^2.$$
\end{prop}
\proof
The freeness of $C$ follows from \cite{DMich,JV}. The fact that
$\tau(C)=\mu(C)=39$ can be seen directly, since $C$ has 3 singularities $A_1$ and 4 ordinary singularities of order 4, which are always quasi-homogeneous, see for instance \cite[Exercise 7.31]{RCS}. The formula for the Alexander polynomial  $\Delta^1_C(t)$ is obtained using any of the two Singular codes described above. From the point of view of pencil-type arrangements developed above, here $\ell=4$ and hence
$\chi(T)=-2$. We apply Theorem \ref{thmPEN} with $\LL=\C_U$ the trivial local system and $\LL'$ the rank one local system on $T$ whose monodromy about each deleted point is multiplication by $\al$, a unity root of order $4$, $\al \ne 1$. It turns out that the inequality in Theorem \ref{thmPEN} is in this case equality.
\endproof
The fact that the Alexander polynomial  $\Delta^1_C(t)$ is non-trivial does not depend on the freeness of the curve $C$, as the following result shows.
\begin{prop}
\label{propT1c} 
Consider the conic-line arrangement 
$$C:f=(y^2-z^2)(x^2-z^2)(x^2-2y^2+z^2)=0$$ 
obtained from a conic pencil of Type I by putting together the 4 lines coming from two of the 3 singular conics in the pencil, plus a smooth member. Then $C$ is a nearly free curve with exponents $(2,4)$ and one has 
$\tau(C)=\mu(C)=18$ and 
$$\Delta^1_C(t)=(t-1)^3(t^3-1).$$
\end{prop}
\proof The nearly freeness is obtained by a direct computation using \verb}Singular}. For the rest,  use the code \verb}mono3_wh()}, but first perform again the substitution $z \mapsto x+y+z$ to make sure the line $L:z=0$ does not contain any singularity of $C$. In this case we have $\ell=3$ and again it turns out that the inequality in Theorem \ref{thmPEN} is  an equality.
\endproof
Now we pass to the case of conic-line arrangements constructed using pencils of Type II. Note that for conic-line arrangements constructed from pencils of type different from Type I, the results in \cite{DMich,JV} can no longer be used to prove the freeness.
\begin{prop}
\label{propT2a} 
Consider the conic-line arrangement 
$$C:f=xy(x-y)(x+y+z)(xy+x^2-y^2+xz-yz)=0$$ 
obtained from a conic pencil of Type II by putting together the 4 lines coming from the 2 singular conics in the pencil, plus a smooth member. Then $C$ is a free curve with exponents $(2,3)$, $\tau(C)=19<\mu(C)=20$, and 
$$\Delta^1_C(t)=(t-1)^3(t^3-1).$$
\end{prop}
\proof
The freeness of $C$ is obtained by a direct computation using \verb}Singular}, which also gives
$\tau(C)=19$. For the remaining claims we use the code \verb}mono3_wh()}. The only non quasi-homogeneous singularity is located at $p=(0:0:1)$.
\endproof
Finally we consider conic-line arrangements constructed using pencils of Type III. The following result can be proved as above.
\begin{prop}
\label{propT3a} 
Consider the conic-line arrangement 
$$C:f=x(x-z)(xy-z^2)(x^2-xz+xy-z^2)=0$$ 
obtained from a conic pencil of Type III by putting together the 2 lines coming from the unique singular conic in the pencil, plus two smooth members. Then $C$ is a free curve with exponents $(2,3)$, 
$\tau(C)=19<\mu(C)=21$, and 
$$\Delta^1_C(t)=(t-1)^2(t^3-1).$$
\end{prop}

\begin{rk}
\label{rkPENCILS} As we have seen,
for all the Alexander polynomials $\Delta^1_C(t)$ of the conic-line arrangements $C$ discussed in this section, the {\it existence} of certain roots $\al$ (but not their {\it exact multiplicity} $m(\al)$) of the Alexander polynomial can be obtained from Corollary \ref{corPEN}. A similar remark applies for most of the conic-line arrangements considered in \cite{JKMS}. These conic-line arrangements are associated to conic pencils of Type IV and V, which contain double lines, hence we have in fact {\it Halphen pencils of index 2}. The reduced pencil-type arrangements $C_{2m}$ in  \cite[Theorem 1.1]{JKMS}, $C_{2m}'$ in  \cite[Theorem 1.3]{JKMS} and $C_{2m}''$ in  \cite[Theorem 1.5]{JKMS} are covered by Theorem \ref{thmPEN2}, while the non-reduced arrangements $C_{2m+1}$ in  \cite[Theorem 1.1]{JKMS} and $C_{2m+1}''$ in  \cite[Theorem 1.5]{JKMS} are covered by Theorem \ref{thmPEN3}.
However, to get the exact formulas for the corresponding Alexander polynomials, one has to use a different approach, and this is done in \cite{JKMS}.
\end{rk}

\section{The family of the Hesse arrangements of conics}

 This construction is described in \cite{Dolgachev, KRS} and we reproduce some of that paper's content here. Let $\lambda \not\in \{1, \omega, \omega^2 \}$ with $\omega^2+\omega + 1=0$. 
 Consider the elliptic curve
 $$E_{\lambda} \, : \, x^{3}+y^{3}+z^{3} - 3\lambda xyz=0,$$
 whose $j$-invariant is given by the formula
 \begin{equation}
\label{eqH0}
j(E_{\lambda})=27 \left( \frac {\lambda(\lambda^3+8)}{\lambda^3-1}\right)^3,
 \end{equation}
 see for instance \cite{PP}. 
 The dual curve $C_{\lambda}$ of the elliptic curve
 $E_{\lambda} $
 is a $9$ cuspidal irreducible sextic, i.e., it has exactly $9$ singularities of type $A_{2}$. Let us denote the set of these $9$ cusps by $\mathcal{P}_{9}(\lambda)$. The coordinates of these points have the following shape:
\begin{equation}
\label{eqH1}
\begin{array}{ccc}

p_{1}(\lambda)=(\lambda:1:1), & p_{4}(\lambda)=(\lambda:\omega:\omega^{2}), & p_{7}(\lambda)=(\lambda:\omega^{2}:\omega)\\
p_{2}(\lambda)=(1:\lambda:1), & p_{5}(\lambda)=(\omega^{2}:\lambda:\omega), & p_{8}(\lambda)=(\omega:\lambda:\omega^{2})\\
p_{3}(\lambda)=(1:1:\lambda), & p_{6}(\lambda)=(\omega:\omega^{2}:\lambda), & p_{9}(\lambda)=(\omega^{2}:\omega:\lambda).
\end{array}
\end{equation}

  It turns out that,  for any fixed $\lambda \not\in  \{-2,-2\omega, -2\omega^2, 0,1, \omega, \omega^2 \}$, that is for all $\lambda \in \C$ such that the curve $E_{\lambda}$ is smooth and $j(E_{\lambda})\ne 0$, the set $\mathcal{P}_{9}(\lambda)$ determines the set of exactly $12$ smooth conics, denoted here by $\mathcal{C}_{12}(\lambda)$, see Proposition \ref{propL}. Each of these 12 conics contains exactly $6$ points from $\mathcal{P}_{9}(\lambda)$. Looking from the perspective of point-conic configurations, we obtain a $(9_{8},12_{6})$-configuration, i.e., through each point from $\mathcal{P}_{9}(\lambda)$ there are $8$ conics from $\mathcal{C}_{12}(\lambda)$ passing through it, and on each conic from $\mathcal{C}_{12}(\lambda)$ there are exactly $6$ points from $\mathcal{P}_{9}(\lambda)$.
 Using any algebra software, we can find explicit equations of these $12$ conics, namely
 \[
\begin{array}{l}
Q_1(\lambda)=C_{1,2,3,4,5,6}(\lambda):\,\,f_1(\lambda)=x^{2}+(\lambda+1)(\omega xy+\omega^{2}xz+yz)+\omega^{2}y^{2}+\omega z^{2}=0,\hfill\\
Q_2(\lambda)=C_{1,2,3,7,8,9}(\lambda):\,\,f_2(\lambda)=x^{2}+(\lambda+1)(\omega^{2}xy+\omega xz+yz)+\omega y^{2}+\omega^{2}z^{2}=0,\hfill\\
Q_3(\lambda)=C_{1,2,4,5,7,8}(\lambda):\,\,f_3(\lambda)=xy-\lambda z^{2}=0,\hfill\\
Q_4(\lambda)=C_{1,2,4,6,8,9}(\lambda):\,\,f_4(\lambda)=x^{2}+(\omega\lambda+1)(xy+\omega xz+\omega yz)+y^{2}+\omega^{2}z^{2}=0,\hfill\\
Q_5(\lambda)=C_{1,2,5,6,7,9}(\lambda):\,\,f_5(\lambda)=x^{2}+(\omega^{2}\lambda+1)(xy+\omega^{2}xz+\omega^{2}yz)+y^{2}+\omega z^{2}=0,\hfill\\
Q_6(\lambda)=C_{1,3,4,5,8,9}(\lambda):\,\,f_6(\lambda)=x^{2}+(\omega\lambda +\omega^{2})(xy+yz+\omega xz)+\omega y^{2}+z^{2}=0,\hfill\\
Q_7(\lambda)=C_{1,3,4,6,7,9}(\lambda):\,\,f_7(\lambda)=-\lambda y^{2}+xz=0,\hfill\\
Q_8(\lambda)=C_{1,3,5,6,7,8}(\lambda):\,\,f_8(\lambda)=x^{2}+(\omega^{2}\lambda +\omega)(xy+yz+\omega^{2}xz)+\omega^{2}y^{2}+z^{2}=0,\hfill\\
Q_9(\lambda)=C_{2,3,4,5,7,9}(\lambda):\,\,f_9(\lambda)=x^{2}+(\lambda +\omega^{2})(xy+xz+\omega^{2}yz)+\omega y^{2}+\omega z^{2}=0,\hfill\\
Q_{10}(\lambda)=C_{2,3,4,6,7,8}(\lambda):\,\,f_{10}(\lambda)=x^{2}+(\lambda +\omega)(xy+xz+\omega yz)+\omega^{2}(y^{2}+z^{2})=0,\hfill\\
Q_{11}(\lambda)=C_{2,3,5,6,8,9}(\lambda):\,\,f_{11}(\lambda)=\lambda x^{2}-yz=0,\hfill\\
Q_{12}(\lambda)=C_{4,5,6,7,8,9}(\lambda):\,\,f_{12}(\lambda)=x^{2}+(\lambda +1)(xy+xz+yz)+y^{2}+z^{2}=0,\hfill
\end{array}
\]
where the indices $i,j,\dots,m$ for the conic $C_{i,j, \ldots ,m}(\lambda)$ mean
that this conic contains the $6$ points $p_{s}(\lambda)$ with $s\in\{i,j, \ldots ,m\}$ from $\mathcal{P}_{9}(\lambda)$.

\begin{prop}
\label{propL} 
Let $\Lambda=\C \setminus  \{-2,-2\omega, -2\omega^2, 0,1, \omega, \omega^2 \}$ and consider the following curve arrangement
$$\mathcal{C}_{12}(\lambda):f(\lambda)=f_1(\lambda) f_2 (\lambda) \cdot \ldots \cdot f_{12}(\lambda)=0.$$ 
Then, for any $\lambda \in \Lambda$, the following two properties hold.

\begin{enumerate}

\item  All the conics $Q_j(\lambda)$ for $j=1, \ldots , 12$ are smooth.

\item  At any point $p_k (\lambda)\in \mathcal{P}_{9}(\lambda)$, any two distinct conics $Q_i(\lambda)$ and $Q_j(\lambda)$ passing through $p_k(\lambda)$ are transversal at $p_k(\lambda)$. In particular, the singularity $(\mathcal{C}_{12}(\lambda), p_k(\lambda))$ is an ordinary singular point of multiplicity 8 and hence 
$$\mu(\mathcal{C}_{12}(\lambda), p_k(\lambda))=49.$$

\end{enumerate}

\end{prop}
\proof
The claims follow by direct computations using, as in our case, \verb}Singular} \cite{Sing}. 
The number of computations can be reduced using Remark \ref{rkHesse0}. The condition $\lambda \not\in \{1, \omega, \omega^2 \}$ is necessary, since we want the curve $E_{\lambda}$ to be smooth and the points in \eqref{eqH1} to be distinct.
Note that for $\lambda=0$, it is clear that each of the conics $Q_3$, $Q_7$ and $Q_{11}$ is not smooth, but union of two lines and the polynomial $f$ has multiple factors. For $\lambda =-2u$, where $u^3=1$, a similar situation occurs, namely the conics $Q_1, Q_2, Q_{12}$ are singular for $u=1$, the conics $Q_5,Q_6,Q_{10}$ are singular for $u=\omega$, and the remaining conics $Q_4,Q_8, Q_9$ are singular for $u=\omega^2$. These special values of $\lambda$ are discussed in the next section.
\endproof

The geometry of the Hesse arrangement of conics is very rich. It turns out that the set  $\Sigma_{12}$ of $12$ double intersection points of $\mathcal{C}_{12}(\lambda)$ determines a set of $9$ lines such that on each line we have exactly $4$ points from $\Sigma_{12}$ -- this is exactly the dual Hesse line arrangement. The corresponding coordinates of points in $\Sigma_{12}$ are
\begin{equation}
\label{eqH2}
\begin{array}{c}
(1:1:1),\,(\omega:\omega^{2}:1),\,(\omega^{2}:\omega:1),\,(1:0:0),\,(0:1:0),\,(\omega^{2}:1:1),\\
(1:\omega^{2}:1),\,(\omega:1:1),\,(1:\omega:1),\,(\omega^{2}:\omega^{2}:1),\,(0:0:1),\,(\omega:\omega:1).
\end{array}
\end{equation}

\begin{rk}
\label{rkHesse0} 
Let $N$ be the cyclic group of order $3$ generated by $n:\PP^2 \to \PP^2$, $n(x:y:z)=(x:\omega y: \omega^2 z)$. Let $H$ be the symmetric group acting on $\PP^2$ by permuting the coordinates $x,y,z$. Then it is easy to check that the semidirect product
$$G= N \rtimes H$$
acts transitively on the set of points $\mathcal{P}_{9}(\lambda)$. As there is at most one conic passing through 6 points, it is clear that the group $G$ permutes the 12 conics in $\mathcal{C}_{12}(\lambda)$ and hence the defining equation $f(\lambda)$ of $\mathcal{C}_{12}(\lambda)$ is invariant under $G$ up to constants in $\C^*$. As an example, the map $n$ takes the set
$\{1,2,3,4,5,6\}$ into the set $\{4,5,6,7,8,9\}$, and hence the conic
$Q_1(\lambda)=C_{1,2,3,4,5,6}(\lambda):f_1(\lambda)=0$ is moved to 
$$Q_{12}(\lambda)=C_{4,5,6,7,8,9}(\lambda):f_1(\lambda)\circ n^{-1}=0.$$
It follows that $f_1(\lambda)\circ n^{-1}=f_{12}(\lambda)$ up to a constant in $\C^*$. As a result, all the properties (e.g. Milnor number, Tjurina number, being an ordinary singularity) of the singularity $(\mathcal{C}_{12}(\lambda),p_j(\lambda))$ do not depend on $j=1, \ldots,12$. The group $G$ acts also on the set $\Sigma_{12}$ of 12 points in \eqref{eqH2}, and there are 4 $G$-orbits, namely
$$ G \cdot (1:1:1), \  G \cdot (1:0:0), \ G \cdot (1:\omega:\omega) \text{ and } G \cdot (1:1:\omega).$$
Hence, to study any singularity $(\mathcal{C}_{12}(\lambda),q)$ for  $q \in \Sigma_{12}$, it is enough to consider just one point $q$ in each of these four $G$-orbits.
\end{rk}

\begin{rk}
\label{rkHesse1} 
The Hesse arrangement of conics
$\mathcal{C}_{12}(\lambda)$ 
for any $\lambda \in \Lambda$, is a free curve with exponents $(7,16)$. 
The singular locus of $\mathcal{C}_{12}$ consists of $9$ ordinary points of multiplicity $8$ given in \eqref{eqH1}, and $12$ double points given in \eqref{eqH2} above. Let us point out here that the points of multiplicity $8$, even if these are ordinary singularities, they are not quasi-homogeneous since again  computations preformed by \verb}Singular} give us that for every $k \in \{1, \ldots ,12\}$ 
 $$45= \tau(\mathcal{C}_{12}(\lambda),p_k(\lambda)) < \mu(\mathcal{C}_{12}(\lambda),p_k(\lambda)) = 49.$$
 For all these results, please consult \cite{PoSz}.
\end{rk}
Recall that the Hesse arrangement of conics
$\mathcal{C}_{12}(\lambda)$ for any  $\lambda \in \Lambda$ is a pencil-type arrangement associated with $\ell=4$ reduced members in a Halphen pencil of index 2, see \cite{Dolgachev}. More precisely, define
$$s_1(\la)= f_1(\la)f_2(\la)f_{12}(\la)=x^6+y^6+z^6+
(\la^3+3\la^2-2)(x^3y^3+y^3z^3+x^3z^3)$$
$$+(-3\la^2-3\la)(x^4yz+xy^4z+xyz^4)+(-3\la^3+9\la+3)x^2y^2z^2,
$$
$$s_2(\la)= f_3(\la)f_7(\la)f_{11}(\la)=-\la^2(x^3y^3+y^3z^3+x^3z^3)+\la(x^4yz+xy^4z+xyz^4)$$
$$+(\la^3-1)x^2y^2z^2,$$
$$s_3(\la)= f_4(\la)f_8(\la)f_{9}(\la)=s_1(\la)+3(1-\omega^2)s_2(\la)$$
and
$$s_4(\la)= f_5(\la)f_6(\la)f_{10}(\la)=s_1(\la)+3(1-\omega)s_2(\la).$$
Then each of the curves $F_j: s_j(\la)=0$ is a union of 3 smooth conics, each two of them meeting in 4 points, of which 3 are in $\PPP_9(\la)$
and one is in $\Sigma_{12}$. Conics from two distinct $F_j$'s meet in 4 points all in $\PPP_9(\la)$. As noted in \cite{Dolgachev}, the 4 sextics $F_1,F_2,F_3$ and $F_4$ are reduced members of the Halphen pencil of index 2 given by
\begin{equation}
\label{eqPEN1}
\PPP(\la): s\cdot s_1(\la)+t\cdot (\la (x^3+y^3+z^3)-(\la^3+2)xyz)^2=0. 
\end{equation}
It follows from \cite[Proposition 5.1]{Dolgachev} that one member of this pencil is a sextic with 9 cusps, and hence the generic member is irreducible.
From now on, when discussing this pencil, we set
$$g_1=s_1(\la) \text{ and } g_2=(\la (x^3+y^3+z^3)-(\la^3+2)xyz)^2.$$
Our Theorem \ref{thmPEN2} (resp. Corollary \ref{corPEN}) implies that  the Alexander polynomial 
$\Delta^1_{\mathcal{C}_{12}(\lambda)}(t)$ admits any primitive root of unity $\al$ of order 8 (resp. any root of unity of order 4) as a root with multiplicity $m(\al) \geq 3$ (resp. $m(\al) \geq 2$). The following result makes this claim more precise.
\begin{thm}
\label{thmHesse} 
The Hesse arrangement of conics 
$\mathcal{C}_{12}(\lambda):f(\lambda)=0,$
for any  $\lambda \in \Lambda$, has the Alexander polynomial given by
$$\Delta^1_{\mathcal{C}_{12}(\lambda)}(t)=(t-1)^{11}(t+1)^2(t^2+1)^2(t^4+1)^3.$$
\end{thm}
\proof
The reader unfamiliar with stratified sets and Thom's Isotopy Lemmas can have a look at the paper \cite{Ra2}, where a proof similar to what follows is given and more details are provided. Alternatively, see \cite[Chapter 1]{D1}.
Let $X=\Lambda \times \PP^2$ and consider the first projection
$${\rm pr}_1:X \to \Lambda.$$
We stratify $X$ using the following connected strata:
\begin{itemize}
\item 21 one-dimensional strata given by
$$X_k=\bigcup_{\lambda \in \Lambda}(\lambda,p_k(\lambda))$$
for $k=1,\ldots,9$ and
$$Y_j=\Lambda  \times p'_j,$$
where $p'_j$ is one of the 12 points in \eqref{eqH2}.
\item 12 two-dimensional strata 
$$Z_j=\bigcup_{\lambda \in \Lambda}\{\lambda\} \times Q_j(\lambda)_0,$$
where $Q_j(\lambda)_0=Q_j(\lambda) \setminus (\mathcal{P}_{9}(\lambda) \cup \Sigma _{12})$, with $\Sigma _{12}$ being the set of 12 points $p'_j$ as in \eqref{eqH2}.
\item a three-dimensional stratum $W$ given by
$$W=\{(\lambda,[x:y:z]) \in X \  | \ f(\lambda)(x,y,z) \ne 0 \}.$$
\end{itemize}

It is easy to check, using Proposition \ref{propL} (2), that for any singular point
 $p(\lambda) \in \mathcal{C}_{12}(\lambda)$ the singularity $(\mathcal{C}_{12}(\lambda),p(\lambda))$ gives rise to a $\mu$-constant family of singularities with $\lambda$ varying in $\Lambda$.
 It is known that a $\mu$-constant family of isolated plane curve singularities is in fact $\mu^*$-constant, see for instance \cite[Remark (1.2.15)]{D1},
and hence it gives rise to a Whitney regular stratification,  \cite[Theorem (1.2.13)]{D1}. It follows that the 2-dimensional strata $Z_j$ are Whitney regular over the 1-dimensional strata, and hence the stratification of $X$ described above is Whitney regular. Since the projection ${\rm pr}_1$ is clearly proper and submersive when restricted to any stratum, we get a
 topologically trivial family via Thom's First Isotopy Lemma, see for instance \cite[Theorem (1.3.5)]{D1}. 
 It follows that the topological type of the pair $(\PP^2,\mathcal{C}_{12}(\lambda))$ does not depend on $\lambda \in \Lambda$, and hence the same applies to the complement
$U(\lambda)=\PP^2 \setminus \mathcal{C}_{12}(\lambda)$. 
Recall the definition of $\al_q$ from \eqref{root}. The multiplicity $m(\al_q)$ of the root $\al_q$ of the Alexander polynomial 
$\Delta^1_{\mathcal{C}_{12}(\lambda)}(t)$
can be computed using the formula \eqref{mult3}, with $\al=\al_q$ and the local system $\LL_{\al}$ denoted by $\LL_{\al}(\la)$ to stress the dependence on $\la$.
It follows that the Alexander polynomial $\Delta^1_{\mathcal{C}_{12}(\lambda)}(t)$
is constant when $\lambda$ varies in $\Lambda$. Indeed, for any two values $\lambda, \lambda' \in \Lambda$, there is a stratified homeomorphism
$$\phi: (\PP^2,\mathcal{C}_{12}(\lambda)) \to (\PP^2,\mathcal{C}_{12}(\lambda')),$$
which implies that the homeomorphism obtained by restriction 
$$\phi: U(\lambda) \to U(\lambda')$$
preserves the rank one local systems used to compute the multiplicity of the eigenvalues of the monodromy $h$, namely 
$$\LL_{\al}(\lambda)=\phi ^{*}(\LL_{\al}(\lambda')).$$
Since the cohomology groups $H^1(U(\lambda),\LL_{\al}(\lambda))$ behave well with respect to pull-back and homotopy equivalence, see
for instance \cite[Remark 2.5.12]{Sheaves}, our claim that the Alexander polynomial $\Delta^1_{\mathcal{C}_{12}(\lambda)}(t)$
is constant follows.

Then we can take, for instance, $\lambda=2$, and the corresponding Alexander polynomial $\Delta^1_{\mathcal{C}_{12}(2)}(t)$ is obtained using the code \verb}monof3(q1,q2)} with $q1=3$ and $q2=24$, since this curve is free by Remark \ref{rkHesse1}.
\endproof

\section{The degenerations of the family of Hesse arrangements of conics}

We have seen in the previous section that there are 7 special values for
$\lambda$ which must be excluded in Theorem \ref{thmHesse}. As explained briefly in the proof of Proposition \ref{propL}, for each of these values of $\lambda$, the defining equation $f(\lambda)=0$ has 3 linear factors which occurs with exponents 2, hence $f(\lambda)=0$ is in fact a multiarrangement $\CC(\lambda)$ formed by 9 smooth conics and 3 double lines.
We denote in this section by $\A(\lambda)$ the non-reduced pencil type conic-line arrangement obtained from $\CC(\lambda)$ by replacing the 3 double lines by the corresponding reduced lines, and by $\B(\lambda)$ the reduced pencil-type subarrangement of $\A(\lambda)$ consisting just of the 9 conics.
These new arrangements $\A(\lambda)$ and $\B(\lambda)$ can be regarded as degenerations of the Hesse arrangement of conics $\mathcal{C}_{12}$ and they have very interesting properties themselves,
as the following results show.

\begin{thm}
\label{thmHesse1} 
The degenerations $\A(\lambda)$ for $\lambda \in  \{-2,-2\omega, -2\omega^2, 0,1, \omega, \omega^2 \}$ of the Hesse arrangement of conics 
are free conic-line arrangements with exponents $(7,13)$.
Moreover, the degenerations $\B(\lambda)$ for $\lambda \in  \{1, \omega, \omega^2 \}$ of the Hesse arrangement of conics 
are free conic arrangements with exponents $(7,10)$. 

\end{thm}
 
\proof
To fix our ideas, consider the case $ \lambda=0$. If we set $ \lambda=0$ in the polynomial $f(\lambda)$, we get
$\CC(0): f(0)=0$, where $ f(0)(x,y,z)=x^2y^2z^2 \cdot g(x,y,z)$ with
$$ g =s_1(0)s_3(0)s_4(0).$$
Indeed, in this case  $s_2(0)=-x^2y^2z^2$.
It follows that one has
   $$\A(0): xyz \cdot g(x,y,z)=0 \text{ and } \B(0):g(x,y,z)=0.$$
 The claim about the freeness of the arrangement $\A(0)$ follows using
a direct computation via \verb}Singular}. One can also check that the arrangement $\B(0)$ is neither free, nor nearly free. The remaining 6 special values of $\lambda$ can be treated exactly in the same way.
\endproof

It turns out that the new conic-line arrangements $\A(\lambda)$ and $\B(\lambda)$ have interesting Alexander polynomials. As a sample of this claim, we have the following.

\begin{thm}
\label{thmHesse2} 
The degeneration $\A(0)$ of the Hesse arrangement of conics 
has 9 ordinary singularities of multiplicity $7$ and $12$ nodes $A_1$. Moreover
$$\tau(\A(0))=309 <336=\mu(\A(0))$$ 
and its Alexander polynomial is given by
$$\Delta^1_{\A(0)}(t)=(t-1)^{9}(t^7-1)^2.$$
The degeneration $\B(0)$ of the Hesse arrangement of conics 
has 9 ordinary singularities of multiplicity $6$ and $9$ nodes $A_1$. Moreover
$$\tau(\B(0))=216 <234=\mu(\B(0))$$ 
and its Alexander polynomial is given by
$$\Delta^1_{\B(0)}(t)=(t-1)^{7}(t^3-1)(t^3+1)^2.$$
In particular, $\Delta^1_{\A(0)}(t)$ (resp. $\Delta^1_{\B(0)}(t))$ has as roots some roots of unity of order $7$ (resp. $6$).
\end{thm}
\proof
The Alexander polynomial for the free conic-line arrangement $\A(0)$ can be computed using the code \verb}monof3(q1,q2)}, with $q1=3$ and $q2=\deg \A(0)=21$. The code \verb}mono3_wh()} takes longer time, but also gives the value $\mu(\A(0))=336$, which can be obtained also by looking at the singularities of $\A(0)$. The 9 ordinary singularities of multiplicity 7 are exactly the points $p_j(0)$ for $j=1,\ldots, 12$ from \eqref{eqH1} and the 12 nodes $A_1$ are exactly the points in \eqref{eqH2}. 
The value $\tau(\A(0))$ comes from the freeness property in Theorem \ref{thmHesse1} using \cite{duP}, or from the code \verb}mono3_wh()}.

All claims relating to the conic arrangement $\B(0)$ can be derived using the code \verb}mono3_wh()}. The $9$ ordinary singularities of multiplicity $6$ are exactly the points $p_j(0)$ for $j = 1, \ldots, 9$ from \eqref{eqH1} and the $9$ nodes $A_1$ are exactly the points in \eqref{eqH2} except the points $(1:0:0), (0:1:0)$ and $(0:0:1)$. 
\endproof

\begin{rk}
\label{rkHesse2P} 
The {\it existence} of roots $\al$ of order 7 and resp. 6 in this case can be seen using our approach with rank one local systems. Consider first the conic arrangement $\B(0)$, which is the union of 3 reduced members in the pencil $\PPP(0)$ from \eqref{eqPEN1}. It follows from \cite[Proposition 5.1]{Dolgachev} that one member of this pencil is a sextic with 9 cusps, and hence the generic member is irreducible.
This pencil has  a double fiber  $F_0$  given by $x^2y^2z^2=0$, hence it is still a Halphen pencil of index 2. To conclude, it is enough to apply
 Theorem \ref{thmPEN2} with $\ell=3$ and $k=2$ and we get $m(\al) \geq 2$.
 For the conic-line arrangement $\A(0)$,  which is non-reduced, we conclude using Theorem \ref{thmPEN3}.
\end{rk}

 \begin{thm}
\label{thmHesse3} 
The degeneration $\A(1)$ of the Hesse arrangement of conics 
has three  singularities of multiplicity 11 and nine nodes $A_1$. Moreover
$$\tau(\A(1))=309 <363=\mu(\A(1))$$ 
and its Alexander polynomial is given by
$$\Delta^1_{\A(1)}(t)=(t-1)^{9}(t^7-1)^2.$$
The degeneration $\B(1)$  of the Hesse arrangement of conics 
has three singularities of multiplicity $9$ and nine nodes $A_1$. Moreover
$$\tau(\B(1))=219 <255=\mu(\B(1))$$ 
and its Alexander polynomial is given by
$$\Delta^1_{\B(1)}(t)=(t-1)^7(t^3-1)(t^3+1)^2.$$
In particular, $\Delta^1_{\A(1)}(t)$ (resp. $\Delta^1_{\B(1)}(t))$ has as roots all the primitive roots of unity of order $7$ (resp. $6$).
\end{thm}
\proof
If we set $\lambda=1$ in the polynomial $f(\lambda)$, since
$$s_1(1)=(x^3+y^3+z^3-3xyz)^2,$$
we get
$$f(1)(x,y,z)=   (x+y+z)^2(x+\omega y+\omega^2 z)^2(x+\omega^2 y+\omega z)^2 \cdot h(x,y,z)$$ with
$$ h= s_2(1)s_3(1)s_4(1).$$
It follows that one has
   $$\A(1):  (x+y+z)(x+\omega y+\omega^2 z)(x+\omega^2 y+\omega z) h(x,y,z)=0 \text{ and } \B(1):h(x,y,z)=0.$$
The new features of this case are the following. First of all, when we set $\lambda=1$, the 9 points in \eqref{eqH1} collapse at 3 points, namely
$$a=(1:1:1)=p_1(1)=p_2(1)=p_3(1),$$
$$b=(\omega:\omega^2:1)=p_4(1)=p_5(1)=p_6(1),$$
and
$$c=(\omega^2:\omega:1)=p_7(1)=p_8(1)=p_9(1),$$
which form a $G$-orbit, and even an $N$-orbit, with the notation from
Remark \ref{rkHesse0}.
Next, three conics degenerate to double lines, namely
$$f_1(1)=\ell_1^2 \text{ with } \ell_1=x+\omega y +\omega^2 z, \  \ f_2(1)=\ell_2^2 \text{ with } \ell_2=x+\omega^2 y +\omega z,$$
and
$$f_{12}(1)=\ell_3^2 \text{ with } \ell_3=x+ y + z.$$
It follows that the 9 conics in $\B(1)$, which are $Q_j(1)$ for $j \notin \{1,2,12\}$, pass all through the points $a,b$ and $c$, which become points of multiplicity 9 for the arrangement $\B(1)$. They are points of multiplicity 11 for the arrangement $\A(1)$, since through each of them pass two of the lines $\ell_1,\ell_2$ and $\ell_3$. The 9 nodes $A_1$ in the arrangements $\A(1)$ and $\B(1)$ are exactly the points in \eqref{eqH2} except the points $a,b,c$. 

All the other claims  follow by using the code
\verb}mono3_wh()}. In particular, the value for $\mu(\A(1))$ (resp. $\mu(\B(1))$), shows that the three points of multiplicity $11$ (resp. $9$) are not ordinary singularities. In fact, their Milnor numbers are
$$\mu(\A(1),a)=\mu(\A(1),b)=\mu(A(1),c)=(\mu(\A(1))-9)/3=118$$ and respectively 
$$\mu(\B(1),a)=\mu(\B(1),b)=\mu(B(1),c)=(\mu(\B(1))-9)/3=82.$$

\endproof

\begin{rk}
\label{rkHesse3P} 
The {\it existence} of roots $\al$ of order 7 and resp. 6 in this case can again be seen using our approach with rank one local systems. First, consider the conic arrangement $\B(1)$, which is the union of 3 reduced members in the pencil $\PPP(1)$ from \eqref{eqPEN1}. Now this pencil has the double fiber $F_0:(x^3+y^3+z^3-3xyz)^2=0$,
hence it is still a Halphen pencil of index 2. 
It is easy to see that the generic member is a sextic curve having exactly 3 ordinary triple points $D_4$, located at $a,b$ and $c$. A simple application of B\'ezout Theorem implies that such a sextic is irreducible. To conclude, it is enough to apply
 Theorem \ref{thmPEN2} with $\ell=3$ and $k=2$ and we get $m(\al) \geq 2$.
For the conic-line arrangement $\A(1)$,  we conclude using Theorem \ref{thmPEN3}.

\end{rk}

\begin{rk}
\label{rkHesse2} 
All the singularities of a line arrangement in $\PP^2$ are both weighted homogeneous  and {\it ordinary multiple points}. In fact, these singularities are {\it homogeneous, which implies that they are ordinary multiple points}.
The conic-line arrangements constructed in \cite{JKMS} have singularities which are {\it not ordinary multiple points}.
On the other hand, the conic-line arrangements constructed from conic pencils of Type I in Propositions \ref{propT1a} and \ref{propT1c} have singularities which are {\it homogeneous, and hence ordinary multiple points}. And the conic-line arrangements in Theorems \ref{thmHesse} and \ref{thmHesse2} have singularities which are
{\it ordinary multiple points, but some of them not weighted homogeneous}.

Finally, the conic arrangement $\B(1)$ and the conic-line arrangement $\A(1)$
have some singularities which are not topologically equivalent to ordinary multiple points. One way to see this using Alexander polynomials is as follows. Note that in the case of the conic arrangement $\B(1)$, which has three singular points of multiplicity $9$,
the Alexander polynomial has roots of order $6$. The local Alexander polynomial of an ordinary point of multiplicity $9$ has only roots of order $9$, as shown in Remark \ref{rk0}. This property is preserved for a singularity which is topologically equivalent  to an ordinary point of multiplicity $9$. Since the other singularities of $\B(1)$ are nodes $A_1$,
with trivial local Alexander polynomial, our claim follows from Theorem \ref{thm0}, (2). A similar approach proves the claim for the conic-line arrangement $\A(1)$, which has three singular points of multiplicity 11.
\end{rk}

\begin{rk}
\label{rkHesse3} 
The conic-line arrangement $\A(1)$ can be called the {\it Cremona--Hesse arrangement} for the following reason. Consider the projective plane $\PP^2$ with coordinates $(u:v:w)$, the Cremona birational morphism
$$\phi: \PP^2 \dashrightarrow \PP^2, \ (u:v:w) \mapsto (vw:uw:uv),$$ and the Hesse line arrangement
\begin{equation} 
\label{He1}
\HH: h(u,v,w)=uvw\left( (u^3+v^3+w^3)^3-27u^3v^3w^3\right)=0.
\end{equation} 
The Hesse line arrangement plays a key role in hyperplane arrangement theory, being the only known 4-net, see \cite{FY}. It has 12 nodes and 9 points of multiplicity 4, it is free with exponents $(4,7)$ and the corresponding Alexander polynomial is
$$\Delta^1_{\HH}(t)=(t-1)^{9}(t^4-1)^2,$$
see \cite{BDS, PS2} or \cite[Theorem 8.19]{DHA}.
Note that
$$h(u,v,w)=uvw(u+v+w)(u+\omega v+\omega ^2 w)(u+\omega ^2 v+ \omega w)$$
$$(\omega u +v+w)(u+\omega v+w)(u+v+\omega w)$$
$$(\omega^2u +v+w)(u+\omega^2 v+w)(u+v+\omega^2 w).$$
Consider the scheme pull-back $\CC'= \phi^{-1}(\HH):h'(u,v,w)=0$,
where
$$h'(u,v,w)=h(vw:uw:uv).$$
Then the scheme $\CC'$ is projectively equivalent to the scheme $\CC(1): f(1)=0$ considered above, that is the two defining ideals map to each other under some automorphism of $\PP^2$. More precisely, one has
\begin{equation} 
\label{He2}
f(1)(x,y,z)=h'(x+\omega y+\omega ^2z, x+\omega ^2y+\omega z, x+y+z).
\end{equation}
It follows that the conic-line arrangement $\A(1)$ is projectively equivalent to the conic-line arrangement $\CC'_{red}$, obtained from the scheme
$\CC'$ by replacing the double lines by reduced lines. In particular, the conic arrangement $\B(1)$ is projectively equivalent to a  conic arrangement obtained by taking the 3 members of the pencil of sextics
$$s(x^3y^3+y^3z^3+x^3z^3)-3tx^2y^2z^2,$$
corresponding to $s=1$ and $t^3=1$. 

\end{rk}

\begin{rk}
\label{rkHesse4} 
Both conic-line arrangements $\A(0)$ and $\A(1)$ are not reduced pencil-type arrangements, i.e., they are not unions of {\it reduced} members of a pencil. Since $\deg \A(0)= \deg \A(1)=21$, there are only two possibilities for  $\A(0)$ or $\A(1)$ to be reduced pencil-type arrangements, namely the following:

\begin{enumerate}

\item  $\A(0)$ or $\A(1)$ are union of 7 reduced members in a pencil of cubics. Since $\A(0)$ or $\A(1)$ have only 3 lines, and no irreducible cubic as components, this case is impossible.

\item  $\A(0)$ or $\A(1)$ are unions of 3 reduced members in a pencil
of curves of degree 7. This would imply, using Theorem \ref{thmPEN}, that the corresponding Alexander polynomials have as roots any cubic root of unity, which is not the case.

\end{enumerate}

\end{rk}

\end{document}